\documentclass{article}
\title{Surjunctivity for cellular automata in Besicovitch spaces}
\author{Silvio Capobianco\footnote{
School of Computer Science, Reykjavik University,
Kringlan 1, 103 Reykjavik, Iceland;
\texttt{silvio@ru.is}}}
\date{}

\usepackage[english]{babel}
\usepackage{amsfonts,amssymb,graphicx}

\newcommand{\Acal}{\ensuremath{\mathcal A}}
\newcommand{\Conf}{\ensuremath{\mathcal{C}}}

\newcommand{\densinf}{\ensuremath{\mathrm{dens\;inf}}}
\newcommand{\denssup}{\ensuremath{\mathrm{dens\;sup}}}
\newcommand{\Edges}{\ensuremath{\mathcal{E}}}

\newcommand{\Neigh}{\ensuremath{\mathcal{N}}}
\newcommand{\Nset}{\ensuremath{\mathbb{N}}}
\newcommand{\Parts}{\ensuremath{\mathcal{P}}}

\newcommand{\Zset}{\ensuremath{\mathbb{Z}}}

\newtheorem{theorem}{Theorem}[section]
\newtheorem{corollary}[theorem]{Corollary}
\newtheorem{definition}[theorem]{Definition}

\newtheorem{lemma}[theorem]{Lemma}
\newtheorem{proposition}[theorem]{Proposition}

\newenvironment{proof}{\emph{Proof.}}{$\Box$ \\}

\begin{document}

\maketitle

\begin{abstract}
The Besicovitch pseudodistance measures the relative size of the set of points
where two functions take different values;
the quotient space modulo the induced equivalence relation
is endowed with a natural metric.
We study the behavior of cellular automata in the new topology
and show that, under suitable additional hypotheses,
they retain certain properties possessed in the usual product topology;
in particular, injectivity implies surjectivity.
\end{abstract}
\textit{Keywords:
cellular automata,
finitely generated groups,
Besicovitch topology,
surjunctivity.}

\section{Introduction}
Cellular automata (CA) are transformations
of the space $\Conf$ of configurations on a grid
that are induced by a finitary rule
applied uniformly to each point of the grid.
Such mappings are characterized by translational invariance
and continuity in the product topology;
however, since such topology makes the shift a \emph{chaotic} map,
no translation invariant distance on $\Conf$ can induce it.

To overcome this problem in the unidimensional case,
Cattaneo, Formenti, Margara and Mazoyer~\cite{cfmm97}
define a pseudodistance on the space $\Conf=\{0,1\}^\Zset$
by taking the sets of the form $U_n=[-n,\ldots,n]$,
and computing, for each $n$,
the upper limit $d_B$ of the \emph{densities}
of the sets of points $x\in U_n$
where two configurations take distinct values.
The quotient space where $c_1$ and $c_2$ are identified
iff $d_B(c_1,c_2)=0$
has topological properties similar to those
possessed by the space of \emph{difference equations}---which,
as pointed out by Toffoli~\cite{tt84},
are a field of application for CA.
Additionally, CA induce transformations on the resulting quotient space
which can in turn provide information on several properties
of the original CA.

In this paper, we apply the ideas from~\cite{cfmm97}
in the much broader context of \emph{finitely generated groups},
where CA can still be defined~\cite{csc06,csms99,ff03};
we do this by linking pseudodistances
to increasing sequences of finite sets
which ultimately cover the whole group.
This is not just for the sake of generality:
we are also trying to shed light on any links
between the properties of $\Zset$ and $\{U_n\}$,
and those of $d_B$.
We then address a question asked in \cite{bfk99}:
is there a connection between surjectivity of CA
and surjectivity of induced map?
Finally, we ask whether \emph{surjunctivity},
i.e., being either surjective or noninjective,
is a property of the induced map as well as of the CA.

A summary of answers to these questions is given in
\begin{theorem} \label{thm_main}
Let $G$ be a finitely generated group of subexponential growth
(e.g., $\Zset^d$);
let $S$ be a finite set of generators for $G$
(e.g., the von Neumann or Moore neighborhood)
and let $U_n\subseteq G$ be the set of reduced words on $S\cup S^{-1}$
having at most length $n$;
let $2\leq|Q|<\infty$,
and let $X$ be the quotient of $Q^G$
with respect to the equivalence relation
$$c_1\sim_B c_2\;\;\mathrm{iff}\;\;
\lim_{n\to\infty}\frac{|\{g\in U_n:c_1(g)\neq c_2(g)\}|}{|U_n|}=0\,,$$
endowed with the topology induced by the distance
$$d_B(x_1,x_2)
=\limsup_{n\in\Nset}\frac{|\{g\in U_n:c_1(g)\neq c_2(g)\}|}{|U_n|}
\;\;,\;\;c_1\in x_1,c_2\in x_2\,.$$
Let $\Acal$ be a cellular automaton over $G$ having set of states $Q$.
\begin{enumerate}
\item[1.] $\Acal$ induces in a natural way
a Lipschitz continuous $F:X\to X$.
\item[2.] $F$ is surjective if and only if $\Acal$ is surjective.
\item[3.] $F$ is injective if and only if it is invertible.
\end{enumerate}
Moreover, if $G$ is of polynomial growth (e.g., $\Zset^d$) then
\begin{enumerate}
\item[4.] $d_B$ is invariant by translations, and
\item[5.] the classes of $\sim_B$ are the same for each choice of $S$.
\end{enumerate}
\end{theorem}

\section{Background} \label{sec_back}
Let $f,g:\Nset\to[0,+\infty)$.
We write $f(n)\preccurlyeq g(n)$
if there exist $n_0\in\Nset$ and $C,\beta>0$
such that $f(n)\leq C\cdot g(\beta n)$ for all $n\geq n_0$;
we write $f(n)\approx g(n)$
if $f(n)\preccurlyeq g(n)$ and $g(n)\preccurlyeq f(n)$.
Observe that, if either $f$ or $g$ is a polynomial,
the choice $\beta=1$ is always allowed.

We indicate the identity of the group $G$ as $1_G$.
Product and inverse are extended to subsets of $G$ elementwise.
If $E\subseteq G$ is finite and nonempty,
the \textbf{closure} and \textbf{boundary} of $X\subseteq G$ w.r.t. $E$
are the sets
$X^{+E}=\{g\in G:gE\cap X\neq\emptyset\}=XE^{-1}$ and 
$\partial_E X=X^{+E}\setminus X$, respectively;
in general, $X\not\subseteq X^{+E}$ unless $1_G\in E$.
$S\subseteq G$ is a \textbf{set of generators}
if the graph $(G,\Edges_S)$,
where $\Edges_S=\{(x,xz):x\in G,z\in S\cup S^{-1}\}$,
is connected;
a group is \textbf{finitely generated} (briefly, f.g.)
if it has a finite set of generators (briefly, f.s.o.g.).
The \textbf{distance} between $g$ and $h$ w.r.t. $S$
is their distance in the graph $(G,\Edges_S)$;
the \textbf{length} of $g\in G$ w.r.t. $S$
is its distance from $1_G$.
The \textbf{disk} of center $g$ and radius $r$ w.r.t. $S$
will be indicated by $D_{r,S}(g)$;
we will omit $g$ if equal to $1_G$,
and $S$ if irrelevant or clear from the context.
Observe that $D_r(g)=gD_r$,
and that $(D_{n,S})^{+D_{R,S}}=D_{n+R,S}$.
For the rest of the paper, we will only consider f.g. infinite groups.

The \textbf{growth function} of $G$ w.r.t. $S$ is $\gamma_S(n)=|D_{n,S}|$.
It is well known~\cite{dlh} that $\gamma_S(n)\approx \gamma_{S'}(n)$
for any two f.s.o.g. $S$, $S'$.
$G$ is of \textbf{subexponential growth}
if $\gamma_S(n)\preccurlyeq\lambda^n$ for all $\lambda>1$;
$G$ is of \textbf{polynomial growth}
if $\gamma_S(n)\approx n^k$ for some $k\in\Nset$.
Observe that, if $G=\Zset^d$, then $\gamma_S(n)\approx n^d$.


A sequence $\{X_n\}\subseteq\Parts(G)$ of \emph{finite} subsets of $G$
is \textbf{exhaustive} if $X_n\subseteq X_{n+1}$ for every $n\in\Nset$
and $\bigcup_{n\in\Nset}X_n=G$.
$\{D_n(g)\}$ is an exhaustive sequence.
For $U\subseteq G$, the \textbf{lower} and \textbf{upper density}
of $U$ w.r.t. the exhaustive sequence $\{X_n\}$ are, respectively,
the lower limit $\densinf_{\{X_n\}}U$
and the upper limit $\denssup_{\{X_n\}}U$
of the quantity $|U\cap X_n|/|X_n|$.
An exhaustive sequence is \textbf{amenable}
or a \textbf{F\o lner sequence}~\cite{dlh,ff03,n62} if
\begin{equation} \label{eq_folner}
\lim_{n\to\infty}\frac{|\partial_E X_n|}{|X_n|}=0
\end{equation}
for every finite $E\subseteq G$;
a group is amenable if it has an amenable sequence.
Observe that 
$\{X_n\}$ is amenable iff it satisfies (\ref{eq_folner})
for all the $E$'s in a single exhaustive sequence,
or for $E=S$ f.s.o.g.
If $G$ is of subexponential growth,
then $\{D_n\}$ contains an amenable subsequence,
and is itself amenable if $G$ is of polynomial growth (cf.~\cite{dlh}).

If $2\leq|Q|<\infty$ and $G$ is a f.g. group,
the space $\Conf=Q^G$ of \textbf{configurations} of $G$ over $Q$
is homeomorphic to the \emph{Cantor set}.
If $E\subseteq G$ is finite,
a \textbf{pattern} over $Q$ with \textbf{support} $E$ is a map $p\in Q^E$.
For $c\in\Conf$, $g\in G$, $c^g\in\Conf$ is defined by
$c^g(h)=c(gh)$ for all $h\in G$;
transformations of $\Conf$ of the form $c\mapsto c^g$ for a fixed $g\in G$
are called \textbf{translations}.
For $G=\Zset$ and $g=+1$,
the translation $c\mapsto c^{+1}$ is the shift map.
A pattern $p\in Q^E$ \textbf{occurs} in $c\in Q^G$
if $(c^g)_{|E}=p$ for some $g\in G$.

A \textbf{cellular automaton} (briefly, CA) over $G$ 
is a triple $\Acal=\left<Q,\Neigh,f\right>$,
where the \textbf{set of states} $Q$ is finite and has at least two elements,
the \textbf{neighborhood index} $\Neigh\subseteq G$ is finite and nonempty,
and the \textbf{local evolution function} $f$ maps $Q^\Neigh$ into $Q$.
The map $F_\Acal:Q^G\to Q^G$ defined by
\begin{equation} \label{eq_ul-def}
(F_\Acal(c))_g=f\left((c^g)_{|\Neigh}\right)
\end{equation}
is the \textbf{global evolution function} of $\Acal$.
Observe that $F_\Acal$ is continuous in the product topology
and commutes with translations.
$\Acal$ is injective, surjective, and so on, if $F_\Acal$ is.

Two patterns $p_1,p_2\in Q^E$ are \textbf{mutually erasable}
(briefly, m.e.) for $\Acal$
if $F_\Acal(c_1)=F_\Acal(c_2)$
for any $c_1,c_2\in\Conf$
such that $(c_1)_{|E}=p_1$, $(c_2)_{|E}=p_2$,
and $(c_1)_{|G\setminus E}=(c_2)_{|G\setminus E}$.
$\Acal$ is \textbf{preinjective} if does \emph{not} have two m.e. patterns.
A pattern $p$ is a \textbf{Garden of Eden} (briefly, GoE) for $\Acal$
if it does \emph{not} occur in $F_\Acal(c)$ for any $c\in\Conf$.
From the compactness of $\Conf$ follows that
a CA has a GoE pattern iff it is nonsurjective.
By \textbf{Moore-Myhill's theorem for amenable groups}~\cite{csms99},
a CA over an amenable group is surjective iff it is preinjective;
in particular, it is \textbf{surjunctive},
i.e., either surjective or noninjective.

A \textbf{pseudodistance} on a set $X$ is a map $d:X\times X\to[0,+\infty)$
satisfying all of the axioms for a distance,
except $d(x,y)>0$ for every $x\neq y$.
If $d$ is a pseudodistance on $X$,
then $x_1\sim x_2\;\;\mathrm{iff}\;\;d(x_1,x_2)=0$
is an equivalence relation,
and $d(\kappa_1,\kappa_2)=d(x_1,x_2)$ with $x_i\in\kappa_i$
is a distance on $X/\sim$.

\section{Besicovitch distances and cellular automata} \label{sec_besicov}

\begin{definition} \label{def_hd}
Let $G$ be a group, let $2\leq|Q|<\infty$,
let $U\subseteq G$ be finite, let $c_1,c_2\in\Conf=Q^G$.
The \emph{Hamming (pseudo)distance between $c_1$ and $c_2$ w.r.t. $U$}
is the quantity $H_U(c_1,c_2)=|\{x\in U: c_1(x)\neq c_2(x)\}|$.
\end{definition}
Observe that $H_U(c_1,c_2)\leq H_{U'}(c_1,c_2)$ if $U\subseteq U'$.
If $U=D_{n,S}$ we write $H_{n,S}(c_1,c_2)$ instead of $H_{D_{n,S}}(c_1,c_2)$.
It is straightforward to prove
\begin{proposition} \label{prop_bd}
Let $\{X_n\}\subseteq\Parts(G)$ be exhaustive.
Then
\begin{equation} \label{eq_bd}
d_{B,\{X_n\}}(c_1,c_2)
=\limsup_{n\in\Nset}\frac{H_{X_n}(c_1,c_2)}{|X_n|}
\end{equation}
is a pseudodistance on $\Conf$,
and is a distance if and only if $G$ is finite.
\end{proposition}
\begin{definition} \label{def_bd}
The quantity (\ref{eq_bd}) is called
the \emph{Besicovitch distance of $c_1$ and $c_2$ w.r.t. $\{X_n\}$}.
The quotient space $\Conf_{B,\{X_n\}}=\Conf/\sim_{B,\{X_n\}}$,
where $c_1\sim_{B,\{X_n\}}c_2$ iff $d_{B,\{X_n\}}(c_1,c_2)=0$,
is called the \emph{Besicovitch space induced by $\{X_n\}$}.
\end{definition}
By an abuse of language, we will also indicate as $\Conf_{B,\{X_n\}}$
the metric space $(\Conf_{B,\{X_n\}},d_{B,\{X_n\}})$.
If $X_n=D_{n,S}$ for some f.s.o.g. $S$,
we write $d_{B,S}$ instead of $d_{B,\{D_{n,S}\}}$,
and speak of \emph{Besicovitch distance w.r.t. $S$};
similar nomenclature and notation shall be used in analogous cases.
Observe that, if $G$ is infinite,
and $c_k(g)=c(g)$ if and only if $g\in X_k$,
then $c_k\to c$ in the product topology,
but $d_{B,\{X_n\}}(c_k,c)=1$ for all $k\in\Nset$,
so that $d_{B,\{X_n\}}$ is not continuous in the product topology.

Definition~\ref{def_bd} is an extension of the one given in~\cite{cfmm97}
for the case $Q=\{0,1\}$, $G=\Zset$, $S=\{+1\}$.
In general, the topology of $\Conf_{B,\{X_n\}}$
is very different from that of $\Conf$:
for example, in the aforementioned case,
$\Conf_{B,S}$ is arcwise connected, not locally compact,
and infinite-dimensional,
while $\Conf$ is totally disconnected, compact, and zero-dimensional.
Also, the equivalence classes of $\sim_{B,\{X_n\}}$
usually depend on $\{X_n\}$;
for example, if $Q=\{0,1\}$, $G=\Zset$,
$X_n=\{-n,\ldots,n\}$, $X'_n=\{-n,\ldots,2^n\}$,
$c_1(g)=0$ for all $g$, $c_2(g)=1$ iff $g<0$,
then $d_{B,\{X_n\}}(c_1,c_2)=1/2$ but $d_{B,\{X'_n\}}(c_1,c_2)=0$.
\begin{theorem} \label{thm_bd-S-inv}
Let $G$ be a group of polynomial growth.
For every $c_1,c_2\in\Conf$,
exactly one of the following happens:
\begin{enumerate}
\item $d_{B,S}(c_1,c_2)>0$ for every f.s.o.g. $S$;
\item $d_{B,S}(c_1,c_2)=0$ for every f.s.o.g. $S$.
\end{enumerate}
\end{theorem}
\begin{proof}
Let $S$ be a f.s.o.g. for $G$
such that $d_{B,S}(c_1,c_2)=0$.
Let $S'$ be another f.s.o.g. for $G$:
there exist $k,n_0\in\Nset$ and $\alpha_1,\alpha_2>0$
such that $\gamma_S(n)\leq\alpha_1n^k$
and $\alpha_2n^k\leq\gamma_{S'}(n)$ for all $n>n_0$.
If $\beta$ satisfies $D_{1,S'}\subseteq D_{\beta,S}$,
then $\gamma_S(\beta n)/\gamma_{S'}(n)\leq\alpha_1\beta^k/\alpha_2$
for all $n>n_0$.
Thus, for all $n$ large enough,
$$\frac{H_{n,S'}(c_1,c_2)}{\gamma_{S'}(n)}
\leq\frac{H_{\beta n,S}(c_1,c_2)}{\gamma_S(\beta n)}
\cdot\frac{\gamma_S(\beta n)}{\gamma_{S'}(n)}
\leq\frac{\alpha_1\beta^k}{\alpha_2}
\cdot\frac{H_{\beta n,S}(c_1,c_2)}{\gamma_S(\beta n)}\;,$$
and the rightmost term vanishes for $n\to\infty$.
\end{proof}
It is proved in~\cite{cfmm97}
that the Besicovitch distance $d_{B,\{+1\}}$ on $\{0,1\}^\Zset$
is invariant by translations;
this is \emph{not} true in the general case.
As a counterexample, let $S=\{a,b\}$
and let $G$ be the free group over $S$;
identify elements of $G$ with reduced words over $S\cup S^{-1}$.
Let $c_1(g)=0$ for all $g\in G$,
and $c_2(g)=1$ if and only if $g$ begins with $a$:
then $c_1^a=c_1$,
but $c_2^a(g)=0$ if and only if $g$ begins with $a^{-1}$,
so that $d_{B,S}(c_1,c_2)=1/4$ but $d_{B,S}(c_1^a,c_2^a)=3/4$.
A generalization of the result in~\cite{cfmm97} is given by
\begin{theorem} \label{thm_bd-amenable-inv}
Let $\{X_n\}\subseteq\Parts(G)$ be such that $\{X_n^{-1}\}$ is amenable.
Then $d_{B,\{X_n\}}$ is invariant by translations.
\end{theorem}
\begin{proof}
Let $S$ be a f.s.o.g. for $G$;
it is sufficient to prove that
$d_{B,\{X_n\}}(c_1^g,c_2^g)=d_{B,\{X_n\}}(c_1,c_2)$
for all $c_1,c_2\in\Conf$, $g\in E=D_{1,S}$.

Given $c\in\Conf$, define $c^-\in\Conf$
as $c^-(g)=c(g^{-1})$ for all $g\in G$.
Then
$$H_{X_n}(c_1^g,c_2^g)=H_{gX_n}(c_1,c_2)
=H_{X_n^{-1}g^{-1}}(c_1^-,c_2^-)
\leq H_{X_n}(c_1,c_2)+|\partial_E X_n^{-1}|$$
for all $g\in E$, $n\in\Nset$,
so that from the amenability of $\{X_n^{-1}\}$
follows $d_{B,\{X_n\}}(c_1^g,c_2^g)\leq d_{B,\{X_n\}}(c_1,c_2)$.
This is true for all $c_1,c_2\in\Conf$, $g\in E$,
so that, by replacing $c_i$ with $c_i^g$ and $g$ with $g^{-1}$,
we get the reverse inequality.
\end{proof}
\begin{corollary} \label{cor_bdinv}
If $\{X_n\}\subseteq\Parts(G)$ is an amenable sequence of symmetric sets,
then $d_{B,\{X_n\}}$ is invariant by translations.
In particular, if $G=\Zset^d$
and $S$ is either the von Neumann or the Moore neighborhood,
then $d_{B,S}$ is invariant by translations.
\end{corollary}
We now ask ourselves which properties do CA possess w.r.t. $d_{B,\{X_n\}}$.
First of all, given $F:\Conf\to\Conf$,
we look after sufficient conditions for
\begin{equation} \label{eq_F-induced}
\overline{F}\left([c]_{\sim_{B,\{X_n\}}}\right)
=\left[F(c)\right]_{\sim_{B,\{X_n\}}}
\end{equation}
to be well defined.
One such condition
is \emph{Lipschitz continuity} w.r.t. $d_{B,\{X_n\}}$,
i.e., existence of $L>0$ such that
\begin{equation} \label{eq_lipsch}
d_{B,\{X_n\}}(F(c_1),F(c_2))\leq L\cdot d_{B,\{X_n\}}(c_1,c_2)
\;\;\forall c_1,c_2\in\Conf\;.
\end{equation}
\begin{theorem} \label{thm_bdeqv}
Let $G$ be a f.g. group
and let $\Acal=\left<Q,\Neigh,f\right>$ be a CA over $G$.
\begin{enumerate}
\item If $\{X_n\}$ is amenable,
then $F_\Acal$ satisfies (\ref{eq_lipsch}) with $L=1+|\Neigh|$.
\item If $\{X_n\}=\{D_{n,S}\}$ for some f.s.o.g. $S$,
and $\Neigh\subseteq D_{r,S}$,
then $F_\Acal$ satisfies (\ref{eq_lipsch}) with $L=(\gamma_S(r))^2$.
\end{enumerate}
\end{theorem}
\begin{proof}
If $X\subseteq G$ and $\Neigh\subseteq E$,
then $H_{X}(F_\Acal(c_1),F_\Acal(c_2))\leq|E|\cdot H_{X^{+E}}(c_1,c_2)$.

If $\{X_n\}$ is amenable, put $E=\Neigh\cup\{1_G\}$: then
$$H_{X_n}(F_\Acal(c_1),F_\Acal(c_2))\leq|E|\cdot H_{X_n^{+E}}(c_1,c_2)
\leq|E|\cdot\left(H_{X_n}(c_1,c_2)+|\partial_E X_n|\right)\;,$$
so that point 1 is achieved because of $\{X_n\}$ being amenable.

If $X_n=D_{n,S}$, put $E=D_{r,S}$: then
$H_{n,S}(F_\Acal(c_1),F_\Acal(c_2))\leq\gamma_S(r)\cdot H_{n+r,S}(c_1,c_2)$,
and since $\gamma_S(n+r)\leq\gamma_S(n)\gamma_S(r)$,
we have for all $n\in\Nset$
$$\frac{H_{n,S}(F_\Acal(c_1),F_\Acal(c_2))}{\gamma_S(n)}
\leq(\gamma_S(r))^2\cdot\frac{H_{n+r,S}(c_1,c_2)}{\gamma_S(n+r)}\;,$$
so that point 2 is achieved by taking upper limits w.r.t. $n$.
\end{proof}
We now define two properties of transformations of $\Conf$
that coincide, respectively, with surjectivity and injectivity
of (\ref{eq_F-induced}), if the latter is well defined.
\begin{definition} \label{def_bd-surj-inj}
Let $\{X_n\}\subseteq\Parts(G)$ be an exhaustive sequence.
$F:\Conf\to\Conf$ is \emph{Besicovitch surjective w.r.t. $\{X_n\}$}
(briefly, \emph{$(B,\{X_n\})$-surjective})
if for all $c\in\Conf$ there exists $c'\in\Conf$ such that
$d_{B,\{X_n\}}(c,F(c'))=0$.
$F$ is \emph{Besicovitch injective w.r.t. $\{X_n\}$}
(briefly, \emph{$(B,\{X_n\})$-injective})
if $d_{B,\{X_n\}}(c_1,c_2)>0$ implies $d_{B,\{X_n\}}(F(c_1),F(c_2))>0$.
\end{definition}
Again, we write $S$ instead of $\{X_n\}$ if the latter is $\{D_{n,S}\}$.
Any surjective function $F$ is also $(B,\{X_n\})$-surjective for all $\{X_n\}$;
however, it is not true \emph{a priori}
that existence of $c'$ such that $d_{B,\{X_n\}}(c,F(c'))=0$
implies existence of $c''$ such that $c=F(c'')$.
As a counterexample, let $\Acal=\left<Q,\Neigh,f\right>$
be a nonsurjective CA over $\Zset$,
let $E\subseteq\Zset$ be finite,
and let $p\in Q^E$ be a GoE pattern for $\Acal$.
Let $k,r\in\Nset$ satisfy
$E\subseteq\{-k,\ldots,k\}$ and $\Neigh\subseteq\{-r,\ldots,r\}$.
Fix $c'\in Q^\Zset$ and
define $c$ by replacing $(F_\Acal(c'))(g)$ with $p(g)$
for all $g\in E$:
then $(F_\Acal(c'))(g)=c(g)$ for all $g\not\in\{-k-r,\ldots,k+r\}$,
so that $d_{B,\{X_n\}}(c,F_\Acal(c'))=0$
for any exhaustive sequence $\{X_n\}$;
however, $c\neq F_\Acal(c'')$ for any $c''\in Q^\Zset$.

It is proved in \cite{bfk99}
that every CA over $\Zset$ with set of states $\{0,1\}$
is surjective if and only if it is $(B,\{+1\})$-surjective.
To extend this result, as we want to do,
we need more tools.
\begin{definition} \label{def_U-W-net}
Let $G$ be a group and let $U,W\subseteq G$ be nonempty.
A \textbf{$(U,W)$-net} is a set $N\subseteq G$
such that the sets $xU$, $x\in N$, are pairwise disjoint,
and $NW=G$.
\end{definition}
Any subgroup is a $(U,U)$-net
for any set $U$ of representatives of its right laterals.
It can be proved via Zorn's lemma~\cite{csc06}
that for every nonempty $U\subseteq G$
there exists a $(U,UU^{-1})$-net;
in particular, for every $R\geq 0$ there exists a $(D_R,D_{2R})$-net.
Observe that any $(U,W)$-net is also a $(Ug,Wh)$-net
for any $g,h\in G$.
Also observe that,
if $N$ is a $(U,W)$-net
and 
$\phi(x)\in xU$ for every $x\in N$,
then $\phi(N)$ is a $(\{1_G\},U^{-1}W)$-net.
\begin{lemma} \label{lemma_E-F-dens}
Let $\{X_n\}$ be amenable
and $N\subseteq G$ be a $(U,W)$-net with $U$ and $W$ finite.
Then $\densinf_{\{X_n\}}\,N\geq 1/|W|$
and $\denssup_{\{X_n\}}\,N\leq 1/|U|$.
\end{lemma}
\begin{proof}
From our observations follows that
it is not restrictive to suppose $1_G\in U\cap W$.
For every $x\in X_n$,
there exist \emph{at most} one pair $(\nu_1,u)\in N\times U$
and \emph{at least} one pair $(\nu_2,w)\in N\times W$
such that $x=\nu_1u$ and/or $x=\nu_2w$:
these imply $\nu_1\in X_nU^{-1}=X_n^{+U}$ and $\nu_2\in X_n^{+W}$,
thus
$$|U|\cdot|N\cap X_n^{+U}|\leq|X_n|\leq|W|\cdot|N\cap X_n^{+W}|\;.$$
If $1_G\in E$, then $X\subseteq X^{+E}$;
hence, $|N\cap X_n|=|N\cap X_n^{+U}|-|N\cap\partial_UX_n|=|N\cap X_n^{+W}|-|N\cap\partial_WX_n|$.
From this and the inequalities above follows
$$\frac{1}{|W|}-\frac{|N\cap\partial_WX_n|}{|X_n|}
\leq\frac{|N\cap X_n|}{|X_n|}
\leq\frac{1}{|U|}-\frac{|N\cap\partial_UX_n|}{|X_n|}\;,$$
and from the amenability of $\{X_n\}$ follows the thesis.
\end{proof}
We are now ready to state and prove the main theorem of this paper.
\begin{theorem} \label{thm_surj-amenable}
Let $\{X_n\}$ be an exhaustive sequence for $G$
that contains an amenable subsequence.
Let $\Acal=\left<Q,\Neigh,f\right>$ be a CA over $G$.
\begin{enumerate}
\item If $\Acal$ is $(B,\{X_n\})$-surjective, then it is surjective.
\item If $\Acal$ is $(B,\{X_n\})$-injective, then it is preinjective.
\item If $\Acal$ is $(B,\{X_n\})$-injective,
then it is $(B,\{X_n\})$-surjective.
\end{enumerate}
\end{theorem}
\begin{proof}
Let $S$ be a f.s.o.g. for $G$.

To prove point 1, suppose, for the sake of contradiction,
that $\Acal$ has a GoE pattern $p$:
it is not restrictive to suppose that the support of $p$
is $D_{k,S}$ for some $k>0$.
Let $N$ be a $(D_k,D_{2k})$-net.
Fix $q\in Q$ and define $c\in\Conf$ as
$$c(g)=\left\{\begin{array}{ll}
p(x^{-1}g)
& \mathrm{if}\;\;g\in D_k(x)\;\;\mathrm{for\;some}\;\;x\in N\;, \\
q & \mathrm{otherwise}\;.
\end{array}\right.$$
Let $c'\in\Conf$.
Let $\phi:N\to G$ be such that, for all $x\in N$,
$\phi(x)\in D_k(x)$ and $(F_\Acal(c'))(\phi(x))\neq c(\phi(x))$:
then $\phi(N)$ is a $(\{1_G\},D_{3k})$-net and
$d_{B,\{X_n\}}(c,F_\Acal(c'))\geq\denssup_{\{X_n\}}\,\phi(N)$.
Let $\{n_j\}$ be such that $\{X_{n_j}\}$ is amenable:
by Lemma~\ref{lemma_E-F-dens},
$$\denssup_{\{X_n\}}\,\phi(N)
\geq\denssup_{\{X_{n_j}\}}\,\phi(N)
\geq\frac{1}{\gamma_S(3k)}\;,$$
so that $d_{B,\{X_n\}}(c,F_\Acal(c'))>0$.
This is true for all $c'\in\Conf$,
therefore $\Acal$ cannot be $(B,\{X_n\})$-surjective.

To prove point 2, suppose, for the sake of contradiction,
that $\Acal$ has two m.e. patterns $p_1,p_2$:
it is not restrictive to suppose that their common support
is $D_{k,S}$ for some $k>0$,
and that $p_1(1_G)\neq p_2(1_G)$.
Let $r\geq 0$ be such that $\Neigh\subseteq D_{r,S}$;
put $R=k+2r+1$.
Let $N$ be a $(D_R,D_{2R})$-net;
fix $q\in Q$ and define $c_1,c_2\in\Conf$ as
$$c_i(g)=\left\{\begin{array}{ll}
p_i(x^{-1}g)
 & \mathrm{if}\;\;g\in D_k(x)\;\;\mathrm{for\;some}\;\;x\in N\;, \\
q\; & \mathrm{otherwise}\;.
\end{array}\right.$$
It is straightforward to check that $F_\Acal(c_1)=F_\Acal(c_2)$.
However, $c_1(x)\neq c_2(x)$ for all $x\in N$:
taking $\{n_j\}$ so that $\{X_{n_j}\}$ is amenable
and reasoning as before, we find
$d_{B,\{X_n\}}(c_1,c_2)\geq\densinf_{\{X_{n_j}\}}\,N\geq 1/\gamma_S(2R)$,
against the hypothesis of $(B,\{X_n\})$-injectivity.

Point 3 follows from point 2,
Moore-Myhill's theorem for amenable groups,
and $(B,\{X_n\})$-surjectivity being implied by surjectivity.
\end{proof}
\begin{corollary} \label{cor_surj-subexp}
Let $G$ be a group of subexponential growth,
let $S$ be a f.s.o.g. for $G$,
and let $\Acal$ be a CA over $G$.
\begin{enumerate}
\item If $\Acal$ is $(B,S)$-surjective,
then it is surjective.
\item If $\Acal$ is $(B,S)$-injective,
then it is preinjective.
\item If $\Acal$ is $(B,S)$-injective,
then it is $(B,S)$-surjective.
\end{enumerate}
\end{corollary}
Observe that, to prove point 3 of Theorem~\ref{thm_surj-amenable},
we do \emph{not} use the fact, implied by Moore-Myhill's theorem,
that injective CA over amenable groups are surjective.
In fact, we deduce preinjectivity directly from $(B,\{X_n\})$-injectivity;
but we do not know (yet) whether this implies injectivity.
At present, our conjecture is that $(B,\{X_n\})$-injectivity
is implied by preinjectivity but does not imply injectivity.
If this were true,
then, for every $\{X_n\}$ containing an amenable subsequence,
any CA would either be
both $(B,\{X_n\})$-injective and $(B,\{X_n\})$-surjective, or neither.

As a final remark, Theorem~\ref{thm_main}
follows from Theorems \ref{thm_bd-S-inv},
\ref{thm_bd-amenable-inv} and \ref{thm_bdeqv},
Corollary~\ref{cor_surj-subexp},
and the observation that surjectivity (injectivity) for $F$
is equivalent to $(B,S)$-surjectivity ($(B,S)$-injectivity) for $\Acal$.

\section{A note on the Weyl distance} \label{sec_weyl}
Given an exhaustive sequence $\{X_n\}$,
we define the \textbf{Weyl (pseudo)distance} as
\begin{equation} \label{eq_wd}
d_{W,\{X_n\}}(c_1,c_2)=\limsup_{n\in\Nset}
\left(\sup_{g\in G}\frac{H_{X_n}(c_1^g,c_2^g)}{|X_n|}\right)\;.
\end{equation}
For $G=\Zset$ and $X_n=[-n,\ldots,n]$,
(\ref{eq_wd}) defines the same quantity as in \cite{bfk99}.
We observe that $d_{W,\{X_n\}}$ is translation invariant whatever $\{X_n\}$ is,
and $d_{W,\{X_n\}}(c_1,c_2)\geq d_{B,\{X_n\}}(c_1,c_2)$
for any two $c_1,c_2\in\Conf$.
However, the metrical properties of $d_W$
are usually worse than those of $d_B$:
in the aforementioned case,
$\Conf_B$ is a complete metric space but $\Conf_W$ is not (cf.~\cite{bfk99}).

We can define $(W,\{X_n\})$-surjectivity
by requiring, for all $c$, the existence of $c'$
such that $d_{W,\{X_n\}}(c,F(c'))=0$.
Then, whatever $\{X_n\}$ is,
$(W,\{X_n\})$-surjectivity implies $(B,\{X_n\})$-surjectivity;
which implies that,
in the hypotheses of Theorem~\ref{thm_surj-amenable},
surjectivity of CA is equivalent to $(W,\{X_n\})$-surjectivity as well.

\section{Acknowledgements}
We thank
Tullio Ceccherini-Silberstein,
Tommaso Toffoli,
Patrizia Mentrasti, and
Henryk Fuk\'s
for the many helpful discussions and encouragements.


\end{document}